# DISCUSSION: ONE-STEP SPARSE ESTIMATES IN NONCONCAVE PENALIZED LIKELIHOOD MODELS[1]

By Cun-Hui Zhang

*Rutgers University*

Penalized methods are commonly used for selecting variables and fitting high-dimensional data. It is well known that the LASSO is biased and thus cannot attain the estimation efficiency of the oracle selector. Recent studies [6, 7, 8, 10, 11] showed that due to the interference of the bias, the LASSO requires quite strong conditions for consistent variable selection. Since the $\ell_1$ penalty has the smallest bias among all convex penalty functions with selection features, these studies naturally draw our attention to methodologies based on concave penalties, or equivalently, nonconcave penalized likelihood.

Frank and Friedman [5] considered the $\ell_\alpha$ penalty for general $\alpha \geq 0$, which is strictly concave for $\alpha < 1$. Their main interest was to use $\alpha$ as a hyperparameter to "bridge" between the subset selection with $\alpha = 0$ and the ridge regression with $\alpha = 2$. Important progresses were made by Fan and Li [3], who advocated the unbiasedness and continuity as essential for variable selectors and carefully developed the SCAD method. In the theoretical front, Fan and Peng [4] proved that the SCAD has the oracle property when the number of variables is a certain fraction power of the sample size. However, nonconcave penalized likelihood methods are still commonly viewed as computationally limited and poorly understood, especially when the number of variables exceeds the number of data points.

Zou and Li made two significant contributions by addressing both the computational and efficiency issues. They developed a fast iterative algorithm for minimizing nonconcave penalized likelihood and proposed a simple one-step method with the oracle property for full rank designs. We congratulate them on this important work. In what follows we relate our work to theirs through discussions on continuity, computational strategies, selection consistency and oracle efficiency.

Received November 2007; revised November 2007.
[1]Supported in part by the NSF Grants DMS-05-04387 and DMS-06-04571 and NSA Grant MDS-904-02-1-0063.

**1. The MC+ method.** Consider penalized squared loss of the form

$$\frac{1}{2n}\|\mathbf{y}-\mathbf{X}\boldsymbol{\beta}\|^2+\sum_{j=1}^{p}\rho(|\beta_j|;\lambda), \tag{1.1}$$

where $\mathbf{y}\in\mathbb{R}^n$ is a response vector, $\mathbf{X}\equiv(\mathbf{x}_1,\dots,\mathbf{x}_p)$ is a design matrix with $p$ covariate vectors $\mathbf{x}_j$, and $\rho(t;\lambda)$ is a penalty function indexed by $\lambda\ge 0$.

In [9] we introduced and studied the MC+, which offers fast, continuous, nearly unbiased and accurate penalized variable selection in high-dimensional linear regression, including the case of $p\gg n$. The MC+ has two elements: a *minimax concave penalty* (MCP) and a *penalized linear unbiased selection* (PLUS) algorithm.

The MCP, given by

$$\rho(t;\lambda)=\lambda\int_0^t\left(1-\frac{x}{\gamma\lambda}\right)^+dx \tag{1.2}$$

with a regularization parameter $\gamma$, minimizes the maximum concavity

$$\kappa(\rho;\lambda)\equiv\max_{t>0}\{-\ddot{\rho}(t;\lambda)\},\qquad \ddot{\rho}(t)\equiv(\partial/\partial t)^2\rho(t;\lambda), \tag{1.3}$$

among all penalty functions satisfying the constraints

$$\dot{\rho}(t;\lambda)=0\qquad\forall t\ge\gamma\lambda,\qquad \dot{\rho}(0+;\lambda)=\lambda, \tag{1.4}$$

where $\dot{\rho}(t;\lambda)\equiv(\partial/\partial t)\rho(t;\lambda)$.

Let $\rho_m(t)$ denote a quadratic spline in $[0,\infty)$ with $m$ knots throughout this discussion, including 0 as a knot. The PLUS computes potentially multiple solutions of the Karush–Kuhn–Tucker-type conditions

$$\begin{cases}\mathbf{x}_j'(\mathbf{y}-\mathbf{X}\widehat{\boldsymbol{\beta}}(\lambda))/n=\operatorname{sgn}(\widehat{\beta}_j(\lambda))\dot{\rho}(|\widehat{\beta}_j(\lambda)|;\lambda), & \widehat{\beta}_j(\lambda)\neq 0,\\ |\mathbf{x}_j'(\mathbf{y}-\mathbf{X}\widehat{\boldsymbol{\beta}}(\lambda))/n|\le\lambda, & \widehat{\beta}_j(\lambda)=0,\end{cases} \tag{1.5}$$

for the possibly nonconvex (1.1), with a penalty of the form $\rho(t;\lambda)=\lambda^2\rho_m(t/\lambda)$. This includes the $\ell_1$ penalty with $m=1$, the MCP with $m=2$ and the SCAD with $m=3$. The output of the PLUS forms a certain *main branch* of the graph of the entire solution set of (1.5). The main branch is a continuous piecewise linear path encompassing from the origin to an "optimal fit" for zero penalty. Other branches of the solution graph form separate loops. The PLUS computes one line segment in the main branch in each step and its computational cost is the same as the LARS per step. For $m=1$, the PLUS becomes the LARS. Moreover, as $\gamma\to\infty$, the MC+ converges to the LASSO for all datasets.

**2. Some simulation results.** We present some simulation results to demonstrate the performance of the LASSO, MC+ and SCAD. Our experiment involves three settings identified by $(n, p, d^o, \beta_*)$, where

$$d^o \equiv \#\{j : \beta_j \neq 0\}, \qquad \beta_* \equiv \min_{\beta_j \neq 0} |\beta_j|, \tag{2.1}$$

characterize the sparsity of the unknown $\boldsymbol{\beta}$. The first two settings are the same as Example 1 of Zou and Li with $\boldsymbol{\beta} = (3, 1.5, 0, 0, 2, 0, 0, 0, 0, 0, 0, 0)'$ and $n = 50$ and 100 respectively. In the third setting, we choose $(n, p, d^o, \beta_*) = (100, 300, 15, 1)$. We divide the 300 components of the $\boldsymbol{\beta}$ vector into 25 continuous blocks, assign the smaller $(3, 1.5, 0, 0, 2, 0, 0, 0, 0, 0, 0, 0)'/1.5$ to 5 random blocks and set entries in other blocks to 0. This makes variable selection harder in the third setting. We generate 1000 independent replications of $(\mathbf{X}, \mathbf{y})$ in each setting, where $\mathbf{y} = \mathbf{X}\boldsymbol{\beta} + \boldsymbol{\varepsilon}$, $\boldsymbol{\varepsilon} \sim N(0, \mathbf{I})$ and $\mathbf{X}$ has i.i.d. Gaussian rows with correlation $2^{-|k-j|}$ between the $j$th and $k$th entries in the same row. For all three procedures, the regularization parameter $\gamma$ is defined as the smallest one for the first part of (1.4) to hold, or equivalently $\gamma = \inf\{t > 0 : \dot{\rho}_m(t) = 0\}$.

We report our results for $\lambda = \sqrt{2(\log p)/n}$ in Table 1. The meaning of MRME is as in Zou and Li. All other entries are averages of the 1000 replications, with ME being the prediction risk as in Zou and Li, MSR being $\|\mathbf{X}\widehat{\boldsymbol{\beta}} - \mathbf{X}\boldsymbol{\beta}\|^2$ as the squared loss for the mean vector, CS being $I\{\widehat{A} = A^o\}$ as the indicator for correct selection, $\mathrm{TM} \equiv |\widehat{A} \setminus A^o| + |A^o \setminus \widehat{A}|$ as the total miss in selection, and the number of PLUS steps $k$ as a measurement of

TABLE 1
*Performance of LASSO, MC+ and SCAD*

| **Method($\gamma$)** | **MRME** | **ME** | **MSE** | **CS** | **TM** | $\boldsymbol{k}$ |
|---|---|---|---|---|---|---|
| | | $n = 50,\ p = 12,\ d^o = 3,\ \beta_* = 1.5$ | | | | |
| LASSO($\infty$) | 0.6129 | 0.2192 | 9.1740 | 0.481 | 0.687 | 4.733 |
| MC+(3.7) | 0.1957 | 0.0753 | 3.3993 | 0.878 | 0.128 | 7.317 |
| SCAD(3.7) | 0.1847 | 0.0689 | 3.1224 | 0.878 | 0.135 | 10.843 |
| | | $n = 100,\ p = 12,\ d^o = 3,\ \beta_* = 1.5$ | | | | |
| LASSO($\infty$) | 0.6794 | 0.1007 | 9.2221 | 0.512 | 0.636 | 4.650 |
| MC+(3.7) | 0.2264 | 0.0361 | 3.4795 | 0.868 | 0.139 | 7.189 |
| SCAD(3.7) | 0.2157 | 0.0327 | 3.1617 | 0.868 | 0.140 | 10.532 |
| | | $n = 100,\ p = 300,\ d^o = 15,\ \beta_* = 1$ | | | | |
| LASSO($\infty$) | | 2.5849 | 148.2765 | 0.000 | 8.271 | 25.227 |
| MC+(2.7) | | 0.2689 | 20.0116 | 0.859 | 0.191 | 41.500 |
| SCAD(3.7) | | 1.3174 | 80.4836 | 0.322 | 1.686 | 59.657 |
| MC+(2.5) | | 0.2373 | 17.9240 | 0.870 | 0.178 | 44.156 |
| SCAD(2.5) | | 0.5510 | 37.6191 | 0.787 | 0.349 | 115.322 |

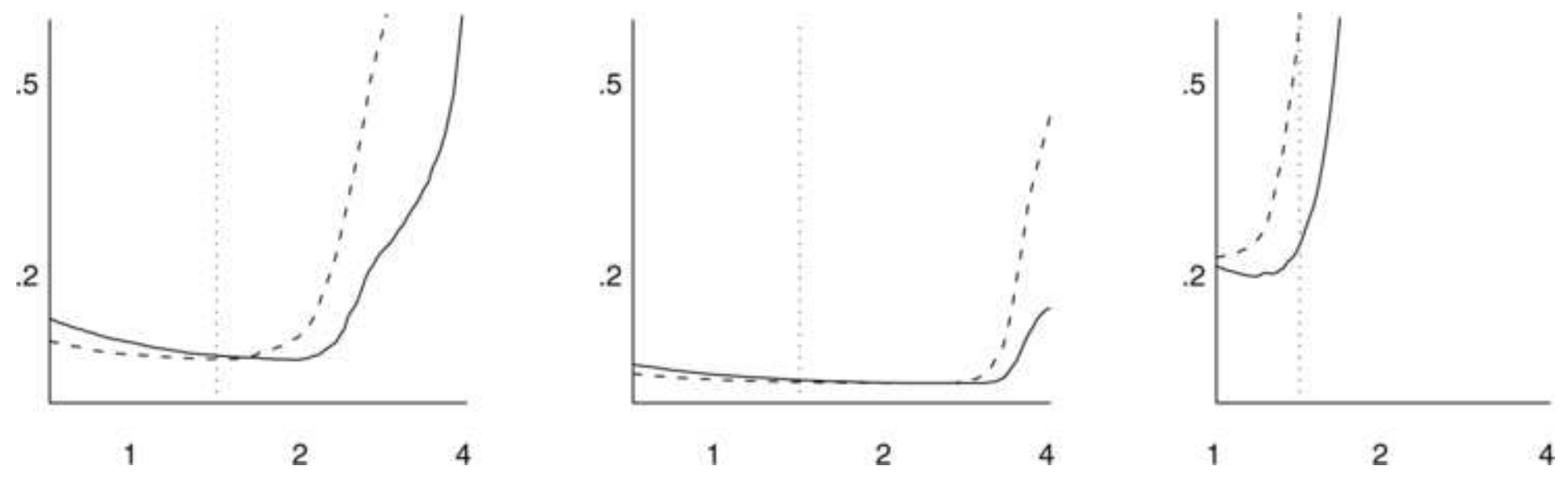

Fig. 1. *The average of 1000 simulated ME of the MC+ (solid) and SCAD (dashed) as functions of* $\lambda/\sqrt{(\log p)/n}$ *for* $(n, p, d^o, \beta_*, \gamma) = (50, 12, 3, 1.5, 3.7)$, $(100, 12, 3, 1.5, 3.7)$ *and* $(100, 300, 15, 1, 2.5)$, *from the left, with dotted vertical at* $\lambda = \sqrt{2(\log p)/n}$.

computational complexity. Here,

$$A^o \equiv \{j : \beta_j \neq 0\} \quad \text{and} \quad \widehat{A} \equiv \{j : \widehat{\beta}_j \neq 0\} \tag{2.2}$$

are the oracle and selected models respectively, and the MRME is undefined for $p > n$. We plot the average of ME and TM against $\lambda$ in Figures 1 and 2. From these results, we observe that the SCAD and MC+ perform similarly in the first two settings, while the MC+ has much stronger performance in the difficult third setting. The LASSO performs poorly in all three settings. This certainly does not represent a thorough simulation comparison of the three methods, but it is consistent with our other simulation experiments [9].

We do not have a definite explanation for the difference in the performance of the SCAD between our simulation and that of Zou and Li in the first two settings, but we offer the following observations. It is clear from Figure 1 that the prediction risk of the SCAD is quite flat in a wide region down to $\lambda/\sqrt{(\log p)/n} = 1$ at least, so that choosing $\lambda$ by CV may cause over fit in view of Figure 2 and the results of Zou and Li. Moreover, 5-fold CV is not designed to choose $\lambda$ accurately unless the dependence of the "optimal" $\lambda$ on $n$ is carefully adjusted, for example, with the factor $n^{-1/2}$. For example, the penalized loss (1.1) with $\rho(t; \lambda) = \lambda|t|$ is equivalent to $\|\mathbf{y} - \mathbf{X}\boldsymbol{\beta}\|^2/2 + \lambda\|\boldsymbol{\beta}\|_1$ for the LASSO with the scale change $\lambda \to n\lambda$, but the best penalty levels in 5-fold CV in these two formulations have effectively 20% difference without adjustment for $n$. Of course, this second problem with CV diminishes if we increase the number of CV folds.

Consider the first two settings in our experiment. In our theorems, a basic upper bound is $2(p - d^o)\Phi(-\lambda\sqrt{n})$ for the probability of selecting some variables with $\beta_j = 0$, which amounts to 0.232 for $\lambda = \sqrt{2(\log p)/n}$. This roughly explains the proportion of incorrect selection 1−CS for the MC+ and SCAD in Table 1. On the other hand, the unbiasedness requires $\beta_* = 1.5 > \gamma\lambda$ with $\gamma = 3.7$, which allows up to $\lambda/\sqrt{(\log p)/n} = 1.82 > \sqrt{2}$ and $2.57 > \sqrt{2}$ respectively for $n = 50$ and 100. Thus, in such cases with a strong signal, the SCAD and MC+ perform better with the larger $\lambda$ as shown in Figures 1 and 2.

**3. Continuity and unbiasedness.** Let $A \subset \{1, \dots, p\}$ represent the model with covariate vectors $\mathbf{x}_j, j \in A$. Define

$$\widehat{\boldsymbol{\beta}}_A(\lambda) \equiv \arg\min \left\{ \frac{1}{2n} \|\mathbf{y} - \mathbf{X}_A \boldsymbol{\beta}_A\|^2 + \sum_{j \in A} \rho(|\beta_j|; \lambda) \right\}, \tag{3.1}$$

where $\boldsymbol{\beta}_A \equiv (\beta_j, j \in A)'$ and $\mathbf{X}_A \equiv (\mathbf{x}_j, j \in A)$. In [9] we prove that, for $\text{rank}(\mathbf{X}_A) = |A|$ and fixed $\lambda$, $\boldsymbol{\beta}_A(\lambda)$ is continuous in $\mathbf{y}$ iff $-\ddot{\rho}(t) < c_{\min}(\mathbf{X}'_A \mathbf{X}_A/n)$ almost everywhere in $t > 0$, where $c_{\min}(\mathbf{M})$ is the smallest eigenvalue of $\mathbf{M}$. Thus, the *global convexity* condition $\kappa(\rho; \lambda) < c_{\min}(\mathbf{X}'\mathbf{X}/n)$ characterizes the global continuity of the global minimizer of (1.1), in the sense of sufficiency and near necessity.

For $p > n$ and sparse $\boldsymbol{\beta}$ with $d^o \equiv \{j : \beta_j \neq 0\} \leq n$, we look for sparse (local) minimizers of (1.1), so that we care about the *sparse continuity* of solutions of (1.5) in the sense of the continuity of (3.1) in $\mathbf{y}$ for all $|A| \leq d^*$, for certain rank $d^* \geq d^o$. This sparse continuity property is characterized by the following *sparse convexity* condition:

$$\kappa(\rho; \lambda) < \min_{|A| = d^*} c_{\min}(\mathbf{X}'_A \mathbf{X}_A/n). \tag{3.2}$$

Under this condition and subject to $\#\{j : \widehat{\beta}_j \neq 0\} \leq d^*/2$, the solution of (1.5) is the unique local minimizer of (1.1) given $\lambda$ and $\mathbf{y}$ and is continuous in $\mathbf{y}$. The global and sparse convexity conditions are not properties of the penalty alone, but they provide the rationale for the use of (1.3) as the measurement of the concavity of the penalty.

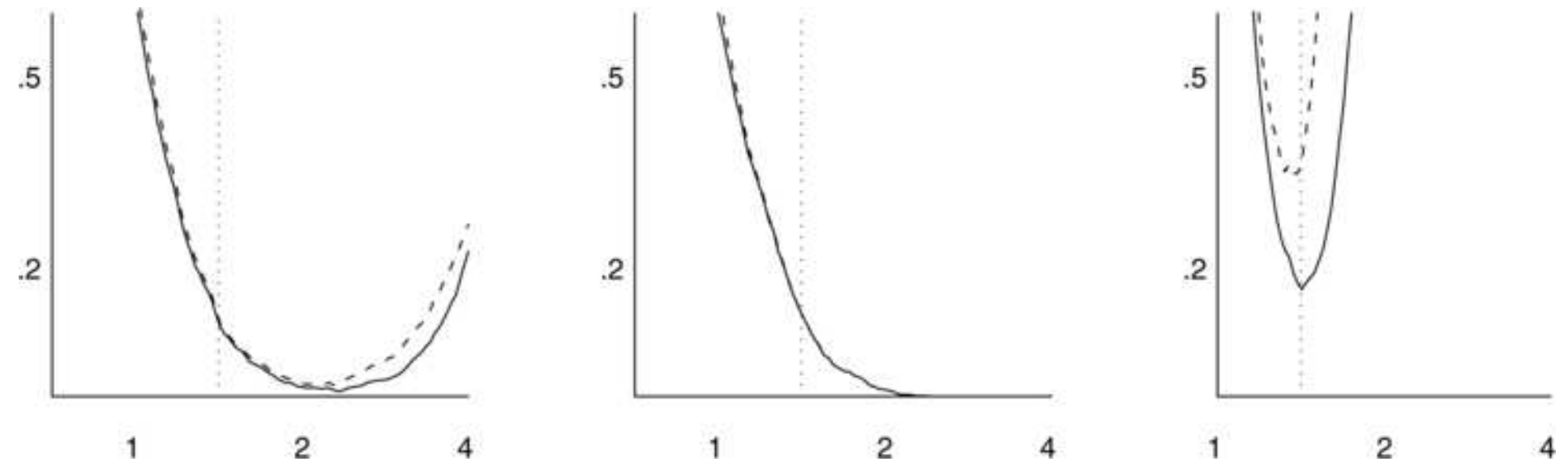

Fig. 2. *The average of 1000 simulated TM of the MC+ (solid) and SCAD (dashed) as functions of $\lambda/\sqrt{(\log p)/n}$ for the experiments in Figure 1.*

The penalty function has selection features if $\dot\rho(0+;\lambda) > 0$. We use the second part of (1.4) to standardize the index $\lambda$ so that it has the interpretation as the threshold for $\beta_j$ for standardized designs with $\|\mathbf{x}_j\|^2/n = 1$. Fan and Li [3] pointed out that penalized estimators are (nearly) unbiased beyond a second threshold $\gamma\lambda$ if the first part of (1.4) holds. Thus, the constraints in (1.4) are natural for unbiased selection. Given (1.4), the MCP provides the sparse continuity for the largest possible rank $d^*$ in (3.2), as it minimizes $\kappa(\rho;\lambda)$. Conversely, given a fixed value of $\kappa(\rho;\lambda)$, the MCP provides the smallest second threshold $\gamma\lambda$ for the unbiasedness, with $\gamma = 1/\kappa(\rho;\lambda)$. Thus, the MCP ensures the continuity and unbiasedness of sparse local minimizers of (1.1) to the greatest extent for general design matrices $\mathbf{X}$. This analysis provide a new point of view since it is clearly different from previous characterizations of penalty functions based on their performance with orthonormal.

**4. The PLUS algorithm.** Consider penalties of the form $\rho(t;\lambda) = \lambda^2\rho_m(t/\lambda)$. Let $\mathbf{z}^* \equiv \mathbf{X}'\mathbf{y}/n$, $\boldsymbol{\chi}_j \equiv \mathbf{x}_j'\mathbf{X}/n$ and $\tau \equiv 1/\lambda$. With the scale change $\mathbf{z} \equiv \tau\mathbf{z}^*$ and $\mathbf{b} \equiv \tau\boldsymbol{\beta}$, (1.5) becomes

$$\begin{cases} z_j - \boldsymbol{\chi}_j'\mathbf{b} = \operatorname{sgn}(b_j)\dot\rho_m(|b_j|), & b_j \neq 0, \\ |z_j' - \boldsymbol{\chi}_j'\mathbf{b}| \le 1 = \dot\rho_m(0+), & b_j = 0. \end{cases} \tag{4.1}$$

The solutions $\mathbf{z} \oplus \mathbf{b}$ of (4.1) form a smooth $p$-dimensional surface $S$ in $\mathbb{R}^{2p}$ composed of $(2m+1)^p$ $p$-parallelepipeds. Given the data $\mathbf{z}^*$, the rescaled solution set of (1.5) is the intersection of this $p$-surface $S$ and the half $p$-subspace $\{(\tau\mathbf{z}^*) \oplus \mathbf{b} : \tau \ge 0\}$. Almost everywhere in $\mathbf{X}$ and $\gamma$, this intersection is composed of a main branch and separate loops. The main branch is piecewise linear, begins with $\mathbf{b} = 0$, and ends with a solution of perfect fit for $\mathbf{z}^*$ [also for $\mathbf{y}$ if $\operatorname{rank}(\mathbf{X}) = n$]. The PLUS algorithm begins with $\mathbf{b} = 0$ and tracks the main branch of the solution set of (1.5) by finding in its $k$th step a second endpoint of its $k$th line segment. Since each step of the PLUS algorithm travels through one distinct parallelepiped in the surface $S$, the algorithm stops in finitely many steps. Our computational strategy for this special nonconvex minimization problem is different from the algorithms of Zuo and Li and other existing iterative ones which converge to a single local minimizer for fixed penalty levels $\lambda$.

Under the global convexity condition, the PLUS finds the unique solution of (1.5) for all $\lambda$ as the global minimizer of (1.1). Otherwise, the value of $\tau = 1/\lambda$ may not be monotone in the PLUS path, so that multiple local minimizers of (1.1) are obtained. We choose the sparsest solution within the PLUS path for a given penalty level. For variable selection purposes, we typically use the *universal penalty level* $\sigma\sqrt{2(\log p)/n}$ or a slightly larger $\lambda$ for large $p$ and standardized designs with $\|\mathbf{x}_j\|^2/n = 1$. We estimate $\sigma$ based on certain *mean residual squares* with a theoretically justified formula for *degrees of freedom*.

Since the PLUS path has to make a turn whenever one of the $b_j$ hits a knot of $\rho_m$, the LASSO is the simplest to compute with $m = 1$, the MC+ is the next for $m = 2$, and then the SCAD for $m = 3$, so on and so forth. The computational complexity is also regulated by the ways the surface $S$ folds as a $p$-vector valued function of $\mathbf{z}$. The orientation of the surface is determined by the eigenvalues of $\mathbf{X}'_A\mathbf{X}_A/n + \mathrm{diag}(\ddot{\rho}(b_j;\lambda), j \in A)$ in individual $p$-parallelepipeds with $A = \{j : b_j \neq 0\}$. Thus, the complexity of the PLUS algorithm is also controlled by the maximum concavity $\kappa(\rho;\lambda)$. For example, when $\kappa(\rho;\lambda) = 1/(\gamma - 1)$ for the SCAD increases from $1/2.7$ to $1/1.5$, the number of required PLUS steps nearly doubles as reported in Table 1.

**5. Selection consistency and oracle estimation efficiency.** A variable selector is consistent if $P\{A^o = \widehat{A}\} \to 1$ with the $A^o$ and $\widehat{A}$ in (2.2). Under selection consistency, efficient estimation after selection yields oracle efficiency under simpler regularity conditions for the $d^o$-dimensional estimation problem as in Theorem 4 of Zou and Li. In this sense, selection consistency implies oracle estimation efficiency.

In [9], we prove that the PLUS solution of (1.5) is selection consistent under mild conditions on $\boldsymbol{\beta}$ and $\mathbf{X}$ in the linear model

$$\mathbf{y} = \mathbf{X}\boldsymbol{\beta} + \boldsymbol{\varepsilon}, \qquad \boldsymbol{\varepsilon} \sim N(0, \sigma^2\mathbf{I}). \tag{5.1}$$

Here we state a simplified version of the theorem. Let $\mathbf{X}_A$ be as in (3.1). The design matrix $\mathbf{X}$ satisfies the sparse Riesz condition if

$$c_* \leq \|\mathbf{X}_A\mathbf{b}\|^2/n \leq c^* \qquad \forall \|\mathbf{b}\| = 1, |A| \leq d^*, \tag{5.2}$$

that is, all the eigenvalues of $\mathbf{X}'_A\mathbf{X}_A/n$ have to lie inside $[c_*, c^*]$ as long as $|A| \leq d^*$. The connection of (5.2) to the Riesz condition on norms was discussed in [10], while sufficient conditions for (5.2) for random matrices were provided in [1, 10], including $d^* = e^{an}$ for fixed positive $\{c_*, c^*, a\}$. The quantities $c_*$ and $c^*$ have been considered as sparse minimum and maximum eigenvalues in [2, 7].

THEOREM 1. *Let* $(\mathbf{X}, \mathbf{y})$ *be as in* (5.1) *with* $\|\mathbf{x}_j\|^2/n = 1$ *and* $\widehat{\boldsymbol{\beta}}^o$ *be the oracle LSE with* $\{j : \widehat{\beta}^o_j \neq 0\} = A^o$ *and* $(\widehat{\beta}^o_j, j \in A^o)' = (\mathbf{X}'_{A^o}\mathbf{X}_{A^o})^{-1}\mathbf{X}'_{A^o}\mathbf{y}$, *where* $A^o$ *is as in* (2.2). *Let* $d^o$ *and* $\beta_*$ *be as in* (2.1). *Let* $\rho(t;\lambda) = \lambda^2\rho_m(t/\lambda)$ *be a quadratic spline penalty function satisfying* (1.4). *Suppose* (5.2) *holds with* $\kappa(\rho_m; 1) < c_*$. *Then, there exist constants* $M_1$ *and* $M_2$ *depending on* $\{c_*, c^*\}$ *and* $\rho_m$ *only, such that for*

$$\beta_* \geq M_1\sigma\sqrt{(1 + \log p)/n}, \qquad M_2 d^o + 1 \leq d^*, \tag{5.3}$$

*and* $\lambda = \sigma\sqrt{2(\log p)/n}$, *the PLUS solution* $\widehat{\boldsymbol\beta} = \widehat{\boldsymbol\beta}(\lambda)$ *of* (1.5) *satisfies*

$$P\{\widehat{A} \neq A^o\} \le P\{\widehat{\boldsymbol\beta} \neq \widehat{\boldsymbol\beta}^o \text{ or } \operatorname{sgn}(\widehat{\boldsymbol\beta}) \neq \operatorname{sgn}(\boldsymbol\beta)\} \to 0 \text{ as } p = p_n \to \infty. \tag{5.4}$$

The selection consistency in Theorem 1 implies that the PLUS solution $\widehat{\boldsymbol\beta}$ achieves the oracle estimation efficiency of $\widehat{\boldsymbol\beta}^o$. Here $\operatorname{sgn}(\boldsymbol\beta)$ means the application of the sign function per component with the convention $\operatorname{sgn}(0) \equiv 0$. We note that $\{p, d^*, d^o, \beta_*\}$ are all allowed to depend on $n$ in Theorem 1. A more general version of Theorem 1 in [9] also allows $(c_*, c^*)$ dependent on $n$, larger $\lambda$ or bounded $p$. An interesting aspect of this result is its validity in cases where $p$ is as large as $e^{an}$ for a fixed small $a > 0$. The one step estimator of Zou and Li can be viewed as adaptive LASSO [12]. As such, it requires an initial estimator which essentially separates the zero and nonzero $\beta_j$ for a certain unspecified threshold.

Department of Statistics and Biostatistics
Hill Center
Busch Campus
Rutgers University
Piscataway, New Jersey 08854
USA
E-mail: czhang@stat.rutgers.edu